\documentclass[12pt]{article}
\usepackage{amsmath, latexsym, amssymb, amscd, float}

\restylefloat{table}
\numberwithin{equation}{section}

\newcommand{\R}{\mathbb R}
\newcommand{\C}{\mathbb C}

\newcommand{\N}{\mathbb N}

\newcommand{\Z}{\mathbb{Z}}

\newcommand{\QED}{\hspace{.2in}\square\newline}
\newcommand{\qed}{\hspace{.2in}\boxminus\newline}

\newtheorem{example}{Example}[section]
\newtheorem{theorem}{Theorem}[section]

\newtheorem{proposition}{Proposition}[section]
\newtheorem{definition}{Definition}[section]
\newtheorem{conjecture}{Conjecture}[section]
\newtheorem{lemma}{Lemma}[section]

\newtheorem{ansatz}{Ansatz}[section]

\begin{document}

\begin{center}
{\Large \textbf{A Gamma Distribution Hypothesis \\ for Prime $k$-tuples}} \vskip 4em

{ J. LaChapelle}\\
\vskip 2em
\end{center}

\begin{abstract}
We conjecture average counting functions for prime $k$-tuples based on a gamma distribution hypothesis for prime powers. The conjecture is closely related to the Hardy-Littlewood conjecture for $k$-tuples but yields better estimates. Possessing average counting functions along with their corresponding exact counting functions allows to implicitly define pertinent $k$-tuple zeta functions. The $k$-tuple zeta functions in turn allow construction of $k$-tuple analogs of explicit formulae. If the zeros of the (implicitly defined) $k$-tuple zeta can be determined, the explicit formulae should yield a (dis)proof of the $k$-tuple analog of the prime number theorem.
\end{abstract}

\vskip 1em

\noindent \emph{Keywords}: Gamma distribution, counting prime $k$-tuples, primes in progressions.

\vskip 1em

\noindent MSC: 11A25, 11N13, 11N60.

\section{Introduction}

The motivation and desire to better understand prime $k$-tuples hardly needs introduction. For a small sample of the literature see \cite{HL,GR,GO,KO,Z} and references therein. Of course there would be no suspense if there existed Euler products for $k$-tuples like that for single primes. The absence of such products seems to indicate that there do not exist generating Dirichlet series whose summands are completely multiplicative in these cases. It is fair to say that this is at the heart of the difficulty in generalizing the single prime case.

This paper hopefully takes a step toward better understanding the distribution of prime $k$-tuples. It is based on two ideas.

The first idea uses the arithmetic function $\mu(n)\Lambda(n)/\log(n)$ to represent the exact prime counting function up to some cut-off $x\in\R_+$,
\begin{equation*}
\pi(x)=-\sum_{n\leq x}\mu(n)\frac{\Lambda(n)}{\log(n)}\;.
\end{equation*}
This simple representation can be readily extended to prime $k$-tuples determined by an admissible set $\mathcal{H}_k=\{0,h_2,\ldots,h_k\}$:
\begin{equation*}
\pi_{(k)}(x)
=(-1)^k\sum_{n\leq x}\mu(n)\cdots\mu(n+h_k)\frac{\Lambda(n)}{\log(n)}
\cdots\frac{\Lambda(n+h_k)}{\log(n+h_k)}\;.
\end{equation*}

The second idea is the gamma distribution hypothesis:  \emph{prime powers} are random variables on the positive reals, and counting them is a random process following a non-homogenous gamma distribution \cite{LA3}. The resulting probability model --- a non-homogenous Poisson process --- yields quite accurate average counting functions associated with the primes.

It is natural to generalize from primes to prime $k$-tuples to test the hypothesis. The obvious tack is to consider a joint gamma distribution on $\R_+^k$. But, in light of the exact $k$-tuple counting function, the counting is modeled by a probability distribution along a ray $\mathbf{R}_{k}\subset\R_+^k$ determined by an admissible set $\mathcal{H}_k=\{0,h_2,\ldots,h_k\}$. Taking this into account leads to an ansatz for the density of \emph{$k$-tuples of prime powers} along $\mathbf{R}_{k}$ up to some cut-off $x$
\begin{equation*}
    {P_{(k)}(n;x)}:=\frac{(-1)^{n-1}}{n!}
    \,\prod_{p< x} \left(1-\frac{\nu_p(\mathcal{H}_k)}{p}\right)\int_{0}^{x}
    (\log(\mathrm{r}))^{n-1}\cdots(\log(\mathrm{r}+h_k))^{n-1}\,d\mathrm{r}
\end{equation*}
where $\nu_p(\mathcal{H}_k)$ is the number of distinct residue classes mod $p$ occupied by the elements in $\mathcal{H}_k$ and the integral is defined by the Cauchy principal value.

There is no reason to expect the probability measure on $\mathbf{R}_{k}$ to coincide with the probability measure on $\R_+$ for the single prime case, and we will argue that
\begin{equation*}
\prod_{p< x} \left(1-\frac{\nu_p(\mathcal{H}_k)}{p}\right)\,d\mathrm{r}=C_{(k)}(x)\,dr
\end{equation*}
where $d\mathrm{r}$ is the measure along $\mathbf{R}_{k}$, $C_{(k)}(x)$ is a normalization defined in (\ref{prime constant}), and $dr$ is the Lebesgue measure on $\R_+$.  For admissible $k$-tuples, the ansatz leads to accurate counting functions because asymptotically $C_{(k)}(x)\sim C_{(k)}$ where $C_{(k)}=\mathfrak{S}(\mathcal{H}_k)$ is the singular series, i.e. the prime $k$-tuple constant.

But the enumeration is secondary. The primary goal is to extract information about prime $k$-tuple distributions, which means we need to discover pertinent $k$-tuple zeta functions implicitly defined by
\begin{equation*}
\log\left(\zeta_{(k)}(s)\right):=\sum_{n=1}^\infty\frac{\Lambda_{(k)}(n)}{\log_{(k)}(n)\,n_{(k)}^s}
\end{equation*}
where $\Lambda_{(k)}(n):=\Lambda(n)\cdots\Lambda(n+h_k)$, similarly $\log_{(k)}(n):=\log(n)\cdots\log(n+h_k)$, and $n_{(k)}$ denotes the geometric mean of the  point $\mathrm{n}_k:=(n, n+h_2,\ldots, n+h_k)\in\R_+^k$. These log-zeta functions are just what one would guess from the structure of the $k$-tuple analogs of the first Chebyshev function defined in the next section.

Unfortunately we haven't found an explicit representation of $\zeta_{(k)}(s)$ that would allow the prime $k$-tuple issue to be settled. But in the final section we motivate the fairly obvious conjecture that
\begin{eqnarray*}
\zeta_{(k)}(s)\,\stackrel{?}{=}\,
\sum_{\mathfrak{n}_{k}}\frac{1}{n_{(k)}^s}
\,\stackrel{?}{=}\,\prod_{p} \left(1-\frac{\nu_p(\mathcal{H}_k)}{p^s}\right) \left(1-\frac{1}{p^s}\right)^{-k}\;.
\end{eqnarray*}
 The sum is over points on an \emph{admissible} ray in the pair-wise coprime $k$-lattice $\mathfrak{N}_+^k\subset\R_+^k$.

This point bears repeating: Possessing both exact counting functions (in terms of standard arithmetic functions) and a model probability distribution facilitates constructing $k$-tuple zeta functions and, subsequently, explicit integral representations of certain counting functions.\footnote{We emphasize that our explicit formulae are left as integral representations. Since we do not determine the complex analytic properties of the $k$-tuple zeta functions, we cannot express the integrals in terms of residues.} To the extent that the $k$-tuple zeta functions and their zeros can be determined, this opens the possibility to attack the problem of  prime $k$-tuple distributions using more-or-less standard methods borrowed from the single prime case.

\section{Counting $k$-tuples}
The first order of business is to collect some arithmetic tools for exact prime $k$-tuple counting.

\begin{proposition}
Let $\mathrm{P}_k$ be the set of prime $k$-tuples, and denote a prime $k$-tuple by $\mathrm{P}_k\ni\mathrm{p}_k=\left(p,p+h_2,\ldots,p+h_k\right)$ with $\mathcal{H}_k:=\{0,h_2,\ldots,h_k\}$ not necessarily admissible. The number of prime $k$-tuples up to some cut-off integer $ x$ is given by\footnote{The subscript $(k)$ is supposed to indicate both the order $k$ of the prime tuple and (implicitly) an associated $\mathcal{H}_k:=\{0,h_2,\ldots,h_k\}$. We will make the dependence on $\mathcal{H}_k$ explicit when necessary.}
\begin{eqnarray}
\pi_{(k)}( x)&:=&\sum_{\mathrm{p}_k\in\mathrm{P}_k;\,p\leq x}1\notag\\
&=&(-1)^k\sum_{n=2}^{ x}\mu(n)\cdots\mu(n+h_k)\frac{\Lambda(n)}{\log(n)}
\cdots\frac{\Lambda(n+h_k)}{\log(n+h_k)}\;.\notag\\
\end{eqnarray}
In particular, the number of prime doubles $(p,p+2i)$ such that $ x-2\geq2i\in\mathbb{N}_+$ is
\begin{equation}
\pi_{(2)}( x):=\sum_{\mathrm{p}_2\in\mathrm{P}_2;\,p\leq x}1=\sum_{n=2}^{ x}\mu(n)\mu(n+2i)\frac{\Lambda(n)}{\log(n)}
\frac{\Lambda(n+2i)}{\log(n+2i)}
\end{equation}
with twin primes corresponding to $i=1$.
\end{proposition}

\emph{Proof}:
Since $\Lambda(n)$ restricts to prime powers $p^\nu$ while $\mu(p^\nu)$ allows only $\nu=1$, then
\begin{equation}
 \mu(n)\Lambda(n)=\left\{\begin{array}{c}
                    -\log(p)\;\;\;\;\;n=p\in\mathrm{P}_1 \\
                    0\;\;\;\;\;\;\;\;\;\;\;\;\mathrm{otherwise}
                  \end{array}\right.\;.
\end{equation}
Loosely, $\mu(n)\Lambda(n)/\log(n)$ acts like a Dirac delta function for primes relative to the discrete measure on the positive reals. More precisely,
\begin{equation}\label{prime sum}
-\sum_{n=2}^{ x}\mu(n)\,\frac{\Lambda(n)}{\log(n)}=\sum_{\mathrm{p}_1\in\mathrm{P}_1;\,p\leq x}1\;.
\end{equation}
Simple induction on $ x$ proves (\ref{prime sum}) since it is obviously true for $ x=2$ and it jumps by one iff $ x+1\in\mathrm{P}_1$.

For the general case let $\mathrm{n}_k:=(n,\ldots,n+h_k)$, then
\begin{eqnarray}
 \mu(n)\Lambda(n)\cdots\mu(n+h_k)\Lambda(n+h_k)=\left\{\begin{array}{l}
(-1)^{k}\log(p)\cdots\log(p+h_k)\,,\;\mathrm{n}_k=\mathrm{p}_k\in\mathrm{P}_{k}\\
0\hspace{2.0in}\mathrm{otherwise}\;.
                         \end{array}\right.\notag\\
\end{eqnarray}
Viewing $\mathrm{n}_k$ as a point in a $k$-lattice and $\mathrm{P}_{k+1}$ as a subset of $\mathrm{P}_{k}\times\mathbb{N}_+=\bigotimes_{k}\mathrm{P}_{1}\times\mathbb{N}_+$, the $k$-tuple result follows after observing that
\begin{eqnarray}
&&\sum_{n=2}^{ x}\left[\frac{\mu(n)\Lambda(n)\cdots\mu(n+h_k)\Lambda(n+h_k)}{\log(n)\cdots\log(n+h_k)}\right]
\frac{\mu(n+h_{k+1})\Lambda(n+h_{k+1})}{\log(n+h_{k+1})}\notag\\
\hspace{.5in}&&=(-1)^k\sum_{\begin{array}{c}
            \scriptstyle{n'\leq x+h_{k+1}} \\
            \scriptstyle{\mathrm{n}_k\in\mathrm{P}_k}
          \end{array}}
\frac{\mu(n')\Lambda(n')}{\log(n')}\,\delta(n'\,,\,(n+h_{k+1}))\notag\\
          \hspace{.5in}&&=(-1)^{k+1}\sum_{\mathrm{p}_{k+1}\in\mathrm{P}_{k+1};\,p\leq x}1\;.
\end{eqnarray}  $\QED$

It is useful to introduce a more compact notation
\begin{equation}
\mu_{(k)}(n):=(-1)^k\mu(n)\cdots\mu(n+h_k)
\end{equation}
and
\begin{eqnarray}
\lambda_{(k)}(n)&:=&\frac{\Lambda(n)\cdots\Lambda(n+h_k)}{\log(n)\cdots\log(n+h_k)}\notag\\
&=:&\frac{\Lambda_{(k)}(n)}{\log_{(k)}(n)}\;.
\end{eqnarray}
So we may write
\begin{equation}
\pi_{(k)}( x)=\sum_{n=2}^{ x}\mu_{(k)}(n)\lambda_{(k)}(n)\;.
\end{equation}

Now define the first and second Chebyshev functions for prime doubles;
\begin{definition}
\begin{eqnarray}
\psi_{(2)}( x)&:=&\frac{1}{2}\sum_{n=2}^{ x}\lambda_{(2)}(n)\log\left(n(n+2i)\right)\;.\\
\theta_{(2)}( x)&:=&\frac{1}{2}\sum_{n=2}^{ x}\mu_{(2)}(n)\lambda_{(2)}(n)\log\left(n(n+2i\right)\;.
\end{eqnarray}
\end{definition}
There are obvious analogs of Chebyshev for higher $k$;
\begin{definition}
\begin{eqnarray}
\psi_{(k)}( x)&:=&\sum_{n=2}^{ x}\lambda_{(k)}(n)\log(n_{(k)})
=\sum_{n=2}^{ x}\frac{\Lambda_{(k)}(n)}{\log_{(k)}(n)}\log(n_{(k)})\\
\theta_{(k)}( x)&:=&\sum_{n=2}^{ x}\mu_{(k)}(n)\lambda_{(k)}(n)\log(n_{(k)})
=\sum_{n=2}^{ x}\mu_{(k)}(n)\frac{\Lambda_{(k)}(n)}{\log_{(k)}(n)}\log(n_{(k)})
\end{eqnarray}
where
\begin{equation}
n_{(k)}:=\left(n(n+h_2)\cdots(n+h_k)\right)^{1/k}\;.
\end{equation}
\end{definition}

\begin{proposition}
\begin{equation}
\theta_{(2)}( x)=\frac{1}{2}\sum_{\mathrm{p}_2\in\mathrm{P}_2;\,p\leq x}\log\left(p(p+2i)\right)
\end{equation}
\end{proposition}

\emph{Proof}: Use the same reasoning as the previous proof. $\QED$

\begin{example}
We can obtain a tight bound on the average (with respect to $i$) prime-double Chebyshev functions:
\begin{eqnarray}
\widehat{\theta_{(2)}}( x)&:=&\frac{\sum_{i=2}^{ x-2}\theta_{(2)}( x)}{\sum_{i=2}^{ x-2}}\notag\\
&=&\frac{1}{2}\frac{1}{( x/2-2)}\sum_{\mathrm{p}_2\in\mathrm{P}_2;\,p\leq x}\left[\log(2^{ x/2-2})+\log(p^{ x/2-2})
+\log\left(\frac{\Gamma(\frac{ x+p}{2})}{\Gamma(\frac{4+n}{2})}\right)\right]\notag\\
&\geq&\frac{1}{2}\sum_{\mathrm{p}_2\in\mathrm{P}_2;\,p\leq x}\left[\log(2)+\log(3)
+\log\left(\frac{\Gamma(\frac{ x+p}{2})}{\Gamma(\frac{4+p}{2})}\right)^{\frac{1}{ x/2-2}}\right]\notag\\
&>&\frac{1}{2}\widehat{\pi_{(2)}}( x)+\frac{\widehat{\pi_{(2)}}( x)}{ x}\left[\log(\Gamma( x/2+1)-1\right]\notag\\
&=&\widehat{\pi_{(2)}}( x)\left[(O(\log( x))+O(1)\right]\;.
\end{eqnarray}
On the other hand,
\begin{eqnarray}
\frac{1}{2}\sum_{\mathrm{p}_2\in\mathrm{P}_2;\,p\leq x}\log\left(p(p+2i)\right)
<\sum_{\mathrm{p}_2\in\mathrm{P}_2;\,p\leq x}\log\left(p+2i\right)
\leq\sum_{\mathrm{p}_2\in\mathrm{P}_2;\,p\leq x}\log( x)
&=&\log( x)\,\pi_{(2)}( x)\;.\notag\\
\end{eqnarray}
So $\widehat{\theta_{(2)}}( x)\asymp\log( x)\,\widehat{\pi_{(2)}}( x)$. Because of Zhang's theorem \cite{Z}, $\widehat{\pi_{(2)}}( x)$ must diverge with $ x$. It follows that $\lim_{ x\rightarrow\infty}\widehat{\theta_{(2)}}( x)/\log( x)=\infty$.  Clearly the same bounds obtain for $\widehat{\psi_{(2)}}( x)$ in terms of $\widehat{J_{(2)}}( x)$ where $J_{(2)}$ is the weighted sum of prime-power doubles.
\end{example}

\section{Prime $k$-tuple conjecture}
This section develops some analytic tools for average $k$-tuple counting based on the gamma distribution hypothesis.

According to \cite{LA3}, events along a directed graph can be modeled by a suitable gamma distribution. For prime $k$-tuples the events occur along a ray $\mathbf{R}_{k}\subset\R_+^k$ determined by an admissible set $\mathcal{H}_k=\{0,h_2,\ldots,h_k\}$.  We are counting prime-power events up to some cut-off $x$, and according to the gamma hypothesis this is a scaled Poisson process.

To learn how to apply the gamma hypothesis for $k$-tuples, let's briefly review the single prime case.  Consider a Poisson process with trivial mean. In this case, counting corresponds to the observation of integers because the events are evenly distributed with unit density. Accordingly, the trivial gamma distribution on $\R_+$ and its associated Poisson process yield a model of the positive integers $\mathbb{Z}_+$ as the cut-off $x\rightarrow\infty$, because they are in one-to-one correspondence with the positive natural numbers $\mathbb{N}_+$.\cite{LA3}

Less trivially, the gamma hypothesis posits that the expected number of prime-power events along $\R_+$ is given by
\begin{equation}
\overline{N(x)}=\sum_{n=1}^\infty\frac{(-1)^{1}}{n!}
    \,\gamma(n,-\log(x))\;.
\end{equation}
 This has the expected form of a scaled Poisson expectation; suggesting we write the incomplete gamma function as an integral in order to infer the associated prime-power probability distribution on $\R_+$:
 \begin{eqnarray}\label{integral}
 \overline{N(x)}&=&\sum_{n=1}^\infty\frac{(-1)^{1}}{n!}\int d\gamma(n,-\log(x))\notag\\
 &=&\sum_{n=1}^\infty\frac{1}{n!}\int(-\log(x))^{n-1}\,dx\notag\\
 &=&-\sum_{n=1}^\infty\frac{(-1)^{n}}{n!}\int_0^x(\log(r))^{n-1}\,dr
 \end{eqnarray}
 where the integral in the third line is taken as the Cauchy principal value. We infer the probability distribution of prime powers on $\R_+$ goes like $\log(r)^{-1}$. Interchanging the sum and integral yields\footnote{We use $\mathrm{Ei}(\log(x))$ instead of $\mathrm{li}(x)$ to remind that the gamma hypothesis applies to the more general case of complex cut-off $x\in\C_+$.} $\overline{N(x)}=\mathrm{Ei}(\log(x))-\log(\log(x))$.

 Return now to the general case. The integral in (\ref{integral}) becomes a multiple integral on $\R_+^k$, and the probability distribution will be a $k$-fold product of logarithms restricted to the appropriate ray determined by $\mathcal{H}_k$. Consequently,  for the general $k$-tuple case  the mean number of prime-power $k$-tuple events is expected to be approximately
\begin{eqnarray}\label{k-density}
\overline{N_{(k)}(x)}&\approx&-\sum_{n=1}^\infty\frac{(-1)^{n}}{n!}\int_{0}^{x}(\log(\mathrm{r}))^{n-1}\cdots (\log(\mathrm{r}+h_k))^{n-1}\,d\mathrm{r}\notag\\
&\approx&-\sum_{n=1}^\infty\frac{(-1)^{n}}{n!}\int_{0}^{x}(\log(r))^{n-1}\cdots (\log(r+h_k))^{n-1}\,dr
\end{eqnarray}
where $\mathrm{r}_k=(\mathrm{r},\mathrm{r}+h_2,\ldots,\mathrm{r}+h_k)\in\mathbf{R}_{k}$ and $r\in\R_+$. In Appendix A we argue these integral representations are equalities if we normalize by
 \begin{equation}
\prod_{p< x} \left(1-\frac{\nu_p(\mathcal{H}_k)}{p}\right)
\end{equation}
and
 \begin{equation}
C_{(k)}(x):=\prod_{p< x} \left(1-\frac{\nu_p(\mathcal{H}_k)}{p}\right)\,\prod_p\left(1-\frac{1}{p}\right)^{-k}
\end{equation}
respectively. This leads to
\begin{ansatz}\label{ansatz}
Let $\mathcal{H}_k=\{0,h_2,\ldots,h_k\}$ be admissible. The expected number of events associated with counting admissible prime-power $k$-tuples up to some cutoff $x$ is
\begin{equation}
    \overline{N_{(k)}(x)}=\sum_{n=1}^\infty{P_{(k)}(n;x)}
:=-\sum_{n=1}^\infty\frac{(-1)^{n}}{n!}
    \,I_{(k)}(n;x)
\end{equation}
where
\begin{eqnarray}
I_{(k)}(n;x)
&=&C_{(k)}(x)\int_{0}^x(\log(r))^{n-1}\cdots(\log(r+h_k))^{n-1}\,dr\notag\\
&=&C_{(k)}(x)\int_{0}^x\log^{n-1}_{(k)}(r)\,dr\notag\\
\end{eqnarray}
with the integral defined by the principal value and $C_{(k)}(x)\sim C_{(k)}$ the singular series.

\end{ansatz}

The analysis in \cite{LA3} suggests that
\begin{equation}
   \overline{N_{(k)}(x)}\approx\sum_{n\leq{x}}\lambda_{(k)}(n)
   -\sum_{n\,\mid\,x}\lambda_{(k)}(n)\;.
 \end{equation}
Recall
\begin{eqnarray}
\lambda_{(k)}(n)&:=&\frac{\Lambda(n)\cdots\Lambda(n+h_k)}{\log(n)\cdots\log(n+h_k)}\notag\\
&=&\frac{\Lambda_{(k)}(n)}{\log_{(k)}(n)}\,\;.
\end{eqnarray}
Moreover, since the sum and integral can be interchanged,
\begin{equation}
    \overline{N_{(k)}(x)}= C_{(k)}(x)\int_{2}^x\frac{1}{\log_{(k)}(r)}\,dr-\mathrm{\;small\; remainder}
\end{equation}
for all $k\in\Z_+$. Again, \cite{LA3} suggests the small remainder term approximates $\sum_{n\,\mid\,x}\lambda_{(k)}(n)$ while
\begin{eqnarray}\label{average J}
  \overline{J_{(k)}(x)}:= C_{(k)}(x)\int_{2}^x\frac{1}{\log_{(k)}(r)}\,dr&=:&C_{(k)}(x)\,\mathrm{Ei}_{(k)}(\log({x}))\notag\\
  &\approx& \sum_{n\leq{x}}\lambda_{(k)}(n)
\end{eqnarray}
is the average $k$-tuple analog of Riemann's counting function.

In particular, for prime doubles
\begin{eqnarray}
    \overline{N_{(2)}(x)}&=&C_{(2)}(x)\sum_{n=1}^\infty\frac{(-1)^{n-1}}{n!}
    \,\int_{2}^x\left(\log(r)\log(r+h_{2i})\right)^{n-1}\,dr\notag\\
&=:&\overline{J_{(2)}(x)}-\overline{\omega_{(2)}(x)}\;.
\end{eqnarray}
where $C_{(2)}(x)$ is the prime-double singular series (which depends on $h_{2i}$).

Given this heuristic motivation, we conjecture:
\begin{conjecture}\label{average primes}
Given an admissible $\mathcal{H}_k=\{0,h_2,\ldots,h_k\}$, the average number of admissible prime $k$-tuples up to some cut-off integer $ x$ is
 \begin{equation}\label{average k-tuples}
 \overline{\pi_{(k)}(x)}=\sum_{m=1}^\infty\frac{\mu(m)}{m}\overline{J_{(k)}(x^{1/m})}
 \end{equation}
 where $\overline{J_{(k)}(x)}:= C_{(k)}(x)\int_{2}^x\frac{1}{\log_{(k)}(r)}\,dr$.
\end{conjecture}
Note that, whereas the Hardy-Littlewood conjecture (\cite{HL} pg. 61, Theorem X) is asymptotic, (\ref{average k-tuples}) holds for all $x>2$. The difference between predicted counts coming from Conjecture \ref{average primes} versus the Hardy-Littlewood conjecture is especially stark for small $x$ or when $x\ll h_k$. This stems from the fact that Hardy-Littlewood is asymptotic and its only $h_i$ dependence comes from the singular series. Appendix B contains some numerical tables illustrating the difference.

Analogous reasoning helps to define the average prime double Chebyshev function:
\begin{definition}
\begin{eqnarray}\label{average Chebyshev}
    \overline{\psi_{(2)}(x)}&:=&C_{(2)}(x)
    \int_{2}^x\frac{\log(r_{(2)})}{\log(r)\log(r+h_{2i})}\,dr\notag\\
    &\approx&\sum_{n\leq x}\frac{\Lambda_{(2)}(n)}{\log_{(2)}(n)}\log(n_{(2)})
\end{eqnarray}
where $r_{(2)}:=\left(r(r+h_{2i})\right)^{1/2}$ is the geometric mean of $(r, r+h_{2i})$.
\end{definition}
This has obvious extensions to higher $k$-tuples:
\begin{definition}
\begin{eqnarray}
    \overline{\psi_{(k)}(x)}&:=&C_{(k)}(x)
    \int_{2}^x\frac{\log(r_{(k)})}{\log_{(k)}(r)}\,dr\notag\\
    &\approx&\sum_{n\leq x}\frac{\Lambda_{(k)}(n)}{\log_{(k)}(n)}\log(n_{(k)})\;.
\end{eqnarray}
\end{definition}

\begin{conjecture}\label{average theta}
\begin{equation}
\overline{\theta_{(2)}}=\sum_{m=1}^\infty\mu(m)\overline{\psi_{(2)}(x^{1/m})}\;.
\end{equation}
\end{conjecture}
Note that $\overline{\psi_{(2)}(x)}
=\tfrac{1}{2}C_{(2)}(x)\left(\mathrm{Ei}(\log({x})+\mathrm{Ei}(\log(x+h))-\mathrm{Ei}(\log(h)\right)
=C_{(2)}(x)\,\mathrm{Ei}(\log({x}))$ follows from (\ref{average Chebyshev}).  Hence, asymptotically, $\overline{\theta_{(2)}(x)}\sim C_{(2)}\left( x/\log( x)\right)$ which is consistent with the Hardy-Littlewood twin prime conjecture.

\section{Explicit formulae}
Having  both exact and average summatory functions allows  to deduce associated $k$-tuple zeta functions and subsequent explicit formulae. Here we will confine attention to prime doubles but indicate the generalization to higher $k$-tuples.

Define the prime-double zeta function \emph{implicitly} by
\begin{definition}
\begin{equation}
\log\left(\zeta_{(2)}(s)\right):=\sum_{n=1}^\infty\frac{\lambda_{(2)}(n)}{n^{s/2}(n+2i)^{s/2}}
=\sum_{n=1}^\infty\frac{\Lambda_{(2)}(n)}{\log_{(2)}(n)\,n_{(2)}^{s}}\;,
\;\;\;\;\Re (s)>1\;.
\end{equation}
\end{definition}
It follows that
\begin{eqnarray}
\log'(\zeta_{(2)}(s))=\frac{\zeta'_{(2)}(s)}{\zeta_{(2)}(s)}
&=&-\sum_{n=1}^\infty\frac{\Lambda_{(2)}(n)}{\log_{(2)}(n)\,n_{(2)}^{s}}\log(n_{(2)})\;.
\end{eqnarray}

Using this log-zeta function, along with the gamma chain  from \cite{LA3} as a guide, we construct an explicit formula  for
  \begin{equation}
\psi_{(2)}(x)=\sum_{n\leq x}\frac{\Lambda_{(2)}(n)}{\log_{(2)}(n)}\log(n_{(2)})\;.
\end{equation}

\begin{proposition}\label{explicit formula}
Put $\widetilde{ x}= x+\epsilon$ with $ x\in\N_+$ and $0<\epsilon<1$. Let $\sigma_a$ be the abscissa of absolute convergence of $\sum_{n=1}^\infty\frac{\lambda_{(2)}(n)\log(n_{(2)})}{n^{s/2}(n+2i)^{s/2}}$. Then, for $c>\sigma_a$,
\begin{eqnarray}\label{explicit eq}
 \psi_{(2)}(x)
 &=&-\lim_{\epsilon\rightarrow0}\lim_{T\rightarrow\infty}\frac{1}{2\pi i}\int_{c-iT}^{c+iT}\Gamma\left(0,-\log(\widetilde{{ \mathsf{x}}}_{(2)}^s)\right)
 \;d\log'(\zeta_{(2)}(s))\,,
 \;\,\;\;\;\;\;c>\sigma_a\notag\\
 &=&\lim_{\epsilon\rightarrow0}\lim_{T\rightarrow\infty}\frac{1}{2\pi i}\int_{c-iT}^{c+iT}\mathrm{Ei}(\log(\widetilde{x}_{(2)}^{s}))
 \;d\log'(\zeta_{(2)}(s))\,,
 \;\,\;\;\;\;\;c>\sigma_a\notag\\
 &=&\sum_{n\leq x}\frac{\Lambda_{(2)}(n)}{\log_{(2)}(n)}\log(n_{(2)})\;.
 \end{eqnarray}
\end{proposition}

\emph{Proof}: First integrate (\ref{explicit eq}) by parts. The boundary term does not contribute because i) a comparison test yields a finite $\sigma_a$ (in fact $\sigma_a=1$) so $\lim_{t\rightarrow\infty}|\log'(\zeta_{(2)}(c+it))|<\infty$ for $c>\sigma_a$;  and ii) $\lim_{t\rightarrow\infty}|\mathrm{Ei}(\log(\widetilde{x}_{(2)}^{s}))|=0$ since
\begin{eqnarray}
\lim_{t\rightarrow\infty}\left|\mathrm{Ei}(\log(\widetilde{x}_{(2)}^{(c+it)})))\right|
&=&\lim_{t\rightarrow\infty}\left|\frac{{ x_{(2)}}^{(c+it)}}{(c+it)\log({ x_{(2)}})}
\left(1+O\left(\frac{1}{(c+it)\log({ x_{(2)}})}\right)\right) \right|\notag\\
&\leq&\frac{{ x_{(2)}}^c}{\log({ x_{(2)}}))}\lim_{t\rightarrow\infty}
\left|\frac{1}{(c+it)}
\left(1+O\left(\frac{1}{(c+it)}\right)\right) \right|=0\;.\notag\\
\end{eqnarray}

Next, following standard arguments, use the truncating integral
\begin{lemma}
\begin{equation}
\frac{1}{2\pi i}\int_{c-i T}^{c+i T}\frac{x^s}{s}ds
=\left\{\begin{array}{l}
1+O\left(\frac{x^c}{T\log(x)}\right)\;\;\;\;x>1 \\
O\left(\frac{x^c}{T\log(x)}\right)\;\;\;\;0<x<1
\end{array}\right.\;.
\end{equation}
\end{lemma}
\emph{proof}: We include the proof for completeness.
For $x>1$ integrate over a rectangle with left edge $(L-i T,L+i T)$ such that $0<L<c$. We have
\begin{equation}
\lim_{L\rightarrow-\infty}\left|\int_{L-i T}^{L+i T}\frac{x^s}{s}ds\right|\leq\lim_{L\rightarrow-\infty}\int_{-T}^{T}\frac{x^L}{|L+ i t|}dt<\lim_{L\rightarrow-\infty}\frac{T x^L}{L}=0\;.
\end{equation}
The top and bottom contribute
\begin{eqnarray}
\left|\int_{-\infty\pm i T}^{c\pm i T}\frac{x^s}{s}ds\right|&\leq&\int_{-\infty}^{0}\frac{-x^{c-r}}{|(c-r)\pm i T|}dr\notag\\
&=&x^c\int_{-\infty}^{0}\frac{-x^{-r}}{|(c-r)\pm i T|}dr\notag\\
&<&x^c\int_{-\infty}^{0}\frac{-x^{-r}}{T}dr\notag\\
&=&\frac{x^{c}}{T\log(x)}\;.
\end{eqnarray}
Finally, the pole at $s=0$ contributes $\mathrm{Res}=1$.

Now, for $x<1$ integrate over the right edge $(R-i T,R+i T)$ with $c< R$. Then
\begin{equation}
\lim_{R\rightarrow\infty}\left|\int_{R-i T}^{R+i T}\frac{x^s}{s}ds\right|\leq\lim_{R\rightarrow\infty}\int_{-T}^{T}\frac{e^{-R|\log(x)|}}{|R+ i t|}dt<\lim_{R\rightarrow\infty}\frac{T e^{-R|\log(x)|}}{R}=0\;.
\end{equation}
The top and bottom contribute the same order as for $x>1$, so the well-known lemma is established. $\qed$

Hence, for $c>\sigma_a$,
\begin{eqnarray}
- \lim_{\epsilon\rightarrow0}\lim_{T\rightarrow\infty}\frac{1}{2\pi i}\int_{c-iT}^{c+iT}\log'(\zeta_{(2)}(s))\frac{\widetilde{{ x}_{(2)}}^{s}} {s}\,ds\notag\\
&&\hspace{-1.25in}=\lim_{\epsilon\rightarrow0}\lim_{T\rightarrow\infty}\frac{1}{2\pi i}\int_{c-iT}^{c+iT}\sum_{n=1}^\infty\frac{\lambda_{(2)}(n)\log(n_{(2)})}{n^{s/2}(n+2i)^{s/2}}
 \frac{\widetilde{{ x}_{(2)}}^{s}}
 {s}\,ds\notag\\
 &&\hspace{-1.25in}=\lim_{\epsilon\rightarrow0}\lim_{T\rightarrow\infty}
 \sum_{n=1}^\infty\frac{\Lambda_{(2)}(n)}{\log_{(2)}(n)}\log(n_{(2)})
 \frac{1}{2\pi i}\int_{c-iT}^{c+iT}
 \frac{\frac{\left(\widetilde{ x}^{1/2}(\widetilde{ x}+2i)^{1/2}\right)^{s}}{n^{s/2}(n+2i)^{s/2}}}
 {s}\,ds\notag\\
 &&\hspace{-1.25in}=\lim_{\epsilon\rightarrow0}\sum_{n\leq\lfloor\widetilde{x}\rfloor}
 \frac{\Lambda_{(2)}(n)}{\log_{(2)}(n)}\log(n_{(2)})\notag\\
 &&\hspace{-1in}=\sum_{n\leq x}\frac{\Lambda_{(2)}(n)}{\log_{(2)}(n)}\log(n_{(2)})
 \end{eqnarray}
where the third equality follows from the lemma. (Justifying the interchange of the sum and integral is straightforward, and interchange of the $T$-limit and sum is allowed because the summand contains $O(n^{-c})$ with $c>1$.)
$\QED$

Clearly this result has teeth only if one possesses an explicit representation of $\zeta_{(2)}(s)$. But if a suitable representation of $\zeta_{(2)}(s)$ can be found and it enjoys analytic properties similar to $\zeta(s)$, then we might expect
\begin{eqnarray}
 \psi_{(2)}(x)&=&\lim_{T\rightarrow\infty}\frac{1}{2\pi i}\int_{c-iT}^{c+iT}
 \mathrm{Ei}(\log({ x}_{(2)}^{s}))\;d\log'(\zeta_{(2)}(s))\notag\\
 &=&\lim_{T\rightarrow\infty}\frac{1}{2\pi i}\int_{c-iT}^{c+iT}\mathrm{Ei}(\log({ x}_{(2)}^{s}))\log''(\zeta_{(2)}(s))\;ds\;.
 \end{eqnarray}
would lead to something like
\begin{equation}
\psi_{(2)}( x)\sim C_{(2)}\mathrm{Ei}(\log(x))
-C_{(2)}\sum_{\rho_{(2)}}\mathrm{Ei}(\log(x^{\rho_{(2)}}))
 +\mathrm{small\,terms}
 \end{equation}
 where the sum would include nontrivial zeros of $\zeta_{(2)}(s)$. Additionally, if non-trivial zeros of $\zeta_{(2)}(s)$ are confined within its critical strip, then the same proof strategy used for the PNT would appear to apply to prime doubles and then all prime $k$-tuples by extension.

 Evidently, if this scenario plays out, then it would seem the prime-double constant $C_{(2)}$ will have to come from the geometric mean $x_{(2)}$ on the coprime lattice and/or $\log''(\zeta_{(2)}(s))$. In this regard, notice that the average $\overline{\psi_{(2)}( x)}$ is given in terms of $C_{(2)}(x)$ and $\int_0^x \log(r_{(2)})\log^{-1}_{(2)}(r)\,dr=\int_0^{x} \log^{-1}(r)\,dr$ whereas the integrand in the explicit formula for $\psi_{(2)}( x)$ contains the factor $\int_0^{x_{(2)}} \log^{-1}(r)\,dr$ (as opposed to $C_{(2)}(x)\int_0^{x} \log^{-1}(r)\,dr$). Figuratively speaking, it's as though $\zeta_{(2)}(s)$ knows how to change $\mathrm{Ei}(\log(x_{(2)}))$ into $C_{(2)}(x)\,\mathrm{Ei}(\log(x))$. This statement is consistent with the interpretation of $C_{(k)}(x)$ as a renormalization required to count prime powers on the pair-wise coprime $k$-lattice (see Appendix A).

Remark that for higher $k$-tuples one should define
\begin{definition}\begin{equation}
\overline{ \varphi_{(k)}(x)}:=C_{(k)}(x)\int_{2}^{x}\frac{\log^{k-1}(r_{(k)})}{\log_{(k)}(r)}\,dr\;.
\end{equation}
\end{definition}
and the $k$-tuple log-zeta function
\begin{definition}
\begin{equation}
\log\left(\zeta_{(k)}(s)\right):=\sum_{n=1}^\infty\frac{\lambda_{(k)}(n)}{n_{(k)}^s}
=\sum_{n=1}^\infty\frac{\Lambda_{(k)}(n)}{\log_{(k)}(n)\,n_{(k)}^s}
\end{equation}
so that
\begin{equation}
\log^{(k-1)'}\left(\zeta_{(k)}(s)\right)
=(-1)^{k-1}\sum_{n=1}^\infty\frac{\Lambda_{(k)}(n)}{\log_{(k)}(n)n_{(k)}^s}\log^{k-1}(n_{(k)})\;.
\end{equation}
\end{definition}

Then to construct an explicit formula at level $k$,  consider
\begin{equation}
\lim_{T\rightarrow\infty}\frac{(-1)^{k-1}}{2\pi i}\int_{c-iT}^{c+iT}\Gamma(0,-\log({ x}_{(k)}^{s}))
 \;d\log^{(k-1)'}\left(\zeta_{(k)}(s)\right)\;.
\end{equation}
Assuming integration by parts to be valid, Perron's formula would apply yielding
\begin{equation}
\varphi_{(k)}(x)=\sum_{n\leq x}\frac{\Lambda_{(k)}(n)}{\log_{(k)}(n)}\log^{k-1}(n_{(k)})\;.
\end{equation}
Finally, assuming favorable analytic properties for $\zeta_{(k)}(s)$ similar to $\zeta(s)$, one would expect to find $\varphi_{(k)}(x)\sim C_{(k)}(x)\mathrm{Ei}(\log(x))$.

\section{Searching for $\zeta_{(k)}(s)$}
The manipulations in the previous section point to a possible representation for $\zeta_{(k)}(s)$ which we formulate as another conjecture. First note that
\begin{eqnarray}
\log\left(\zeta_{(2)}(s)\right)
&>&\sum_{\mathrm{p}_2\in\mathrm{P}_2}
          \sum_{\omega,\omega':\omega=\omega'}\frac{1}{\omega p^{\omega s/2}}\frac{1}{\omega' (p^\omega+2i)^{\omega's/2}}
          \;,\;\;\;\;\;\;\Re(s)>1\notag\\\notag\\
&>&\sum_{\mathrm{p}_2\in\mathrm{P}_2}
          \sum_{\omega,\omega':\omega=\omega'}\frac{1}{\omega (p^{\omega'})^{\omega s/2}}\frac{1}{\omega' ((p+2i)^\omega)^{\omega's/2}}\notag\\\notag\\
&=&\sum_{\mathrm{p}_2\in\mathrm{P}_2}\sum_{\omega^2}\frac{1}{\omega^2 p_{(2)}^{s\,\omega^2}}\notag\\
&=&-\sum_{\mathrm{p}_2\in\mathrm{P}_2}\log\left(1-p_{(2)}^{-s}\right)\;.
\end{eqnarray}
This suggests to define $Z_{(k)}(s):=\prod_{\mathrm{p}_k\in\mathrm{P}_k}\left(1-p_{(k)}^{-s}\right)^{-k}$ which (together with the appendix) motivates
\begin{conjecture}
Let $\mathfrak{R}_k$ be an \textbf{admissible} ray in the pair-wise coprime $k$-lattice $\mathfrak{N}_+^k$ and $n_{(k)}$ the geometric mean of a point $\mathfrak{n}_k=(n, n+h_2, \ldots, n+h_k)\in\mathfrak{R}_k$. Then
\begin{eqnarray}
\zeta_{(k)}(s)
\stackrel{?}{=}\sum_{\mathfrak{n}_{k}}\frac{1}{n_{(k)}^s}
\stackrel{?}{=}\prod_{\mathrm{p}_k\in\mathrm{P}_k}\left(1-p_{(k)}^{-s}\right)^{-k}
\stackrel{?}{=}\prod_{p} \left(1-\nu_p(\mathcal{H}_k)\,p^{-s}\right)\left(1-p^{-s}\right)^{-k}
\end{eqnarray}
where the sum is over all points along an admissible ray $\mathfrak{R}_k\subset\mathfrak{N}_+^k$.
\end{conjecture}

If the conjecture is correct, the $k$-tuple zeta function appears to be a restriction of the multiple zeta function which motivates
\begin{conjecture}
$\zeta_{(k)}(s)$ is meromorphic on $\C$, and the singular part of $(-1)^k\log^{(k)'}(\zeta_{(k)}(s))$ is given by
\begin{equation}
\frac{1}{s-1}\,\frac{(-1)^k}{2\pi i}\oint\log^{(k)'}(\zeta_{(k)}(s))\,ds=\frac{C_{(k)}}{s-1}\;.
\end{equation}
\end{conjecture}

This conjecture is equivalent to the Hardy-Littlewood prime $k$-tuple conjecture in the following sense. If we believe the gamma hypothesis, then being a sum of $1/n_{(k)}^s$ along a ray in the pair-wise coprime $k$-lattice \emph{strongly suggests} $\zeta_{(k)}(s)$ has a first order pole at $s=1$ and there are no other poles, while its zeros are determined by conspiring projections of almost periodic exponentials. Then, together with the explicit formula at level $k$, the conjecture implies
\begin{eqnarray}
{\varphi_{(k)}(x)}
 &=&\left.\lim_{T\rightarrow\infty}\frac{(-1)^k}{2\pi i}\int_{c-iT}^{c+iT}
 \mathrm{Ei}(\log({ x}_{(k)}^{s}))\log^{(k)'}(\zeta_{(k)}(s))\;ds\right|_{s=1}\notag\\\notag\\
 &\sim&C_{(k)}(x)\,\mathrm{Ei}(\log(x))\;.
 \end{eqnarray}
Equivalently,
\begin{equation}
\lim_{N\rightarrow\infty}\frac{1}{N}\sum_{n\leq N}\frac{\Lambda_{(k)}(n)}{\log_{(k)}(n)}\log^k(n_{(k)})
= C_{(k)}\;.
\end{equation}

Of course, possessing poles and zeros of $\zeta_{(k)}(s)$ would be tantamount to evaluating the exact summatory functions. And their singular part would presumably furnish the prime $k$-tuple constants. The goal would be to express the integral in Proposition \ref{explicit formula} as a sum over $k$-tuple zeta residues in the usual way; which would presumably verify Conjecture \ref{average primes} and validate the gamma hypothesis.

\section{Conclusion}
To conclude, it is appropriate to draw attention back to the foundational approach of the pair-wise coprime $k$-lattice. The perspective it affords has i) delivered exact arithmetic and accurate average counting functions for prime $k$-tuples, ii) suggested the Riemann zeta function is but one in a family of $k$-tuple zeta functions, and iii) offered an elementary albeit heuristic explanation for the singular series (see Appendix A).

By its nature, the pair-wise coprime lattice incorporates some simple sieving. This sieving (along with the gamma distribution hypothesis for prime powers) is at the heart of the results presented here; which were derived by elementary means. One can anticipate that valuable results might be gleaned by first applying sophisticated sieve techniques along rays in the pair-wise coprime lattice and then projecting onto $\R_+$ --- as opposed to the other way around.

\vspace{.5in}\noindent
\textbf{Acknowledgments:}
I thank A. Granville for offering some helpful advice at an early stage of this work.

\appendix
\section{The $k$-tuple normalization}

In this appendix we argue that $C_{(k)}(x)$ is asymptotically the prime $k$-tuple constant. Since the argument is based on the gamma distribution hypothesis, it is very close to the standard probability argument. However, the details are a bit different since the counting occurs in the pair-wise coprime lattice.

According to the probability model \cite{LA3} for the case $k=1$, there is a distinction between the positive integer lattice $\mathbb{Z}_+$ and the strictly positive natural numbers $\mathbb{N}_+$ that characterize the counting process associated with the trivial gamma distribution on $\R_+$. However, the trivial gamma distribution yields a Poisson process with unit integer density, and this allows for a meaningful identification between $\mathbb{Z}_+$ and $\mathbb{N}_+$.

But in general the (constrained) probability model represents a prime-power counting process along a ray $\mathbf{R}_{k}$ in the pair-wise coprime $k$-lattice, and there is no guarantee that the trivial gamma distribution along $\mathbf{R}_{k}$ leads to a unit integer density. Accordingly, as a result of restricting to the pair-wise coprime $k$-lattice, it may be necessary to renormalize the probability distribution along $\mathbf{R}_{k}$ if one wants to compare counting processes and maintain the natural identification between $\mathbb{N}_+$ and $\mathbb{Z}_+$.

In other words, comparing exact counting functions to average counting functions in the single prime case, we can interpret
\begin{equation}
\sum_{n\leq x}\frac{\Lambda(n)}{\log(n)}\stackrel{\propto}{\longrightarrow}\int_0^x\frac{dr}{\log(r)}
\end{equation}
as a representation of the averaging process; and then $dr$ is the integrator along $\R_+$ associated with the probability measure of prime powers --- a kind of smoothed $\Lambda(n)$. For the general case we have
\begin{equation}
\sum_{n\leq x}\frac{\Lambda_{(k)}(n)}{\log_{(k)}(n)}\stackrel{\propto}{\longrightarrow}
\int_{0}^{x}\frac{d\mathrm{r}}{\log_{(k)}(\mathrm{r})}\;.
\end{equation}
But now $d\mathrm{r}$ encodes both the `smoothed' $\Lambda_{(k)}(n)$ and the density of prime powers along $\mathbf{R}_{k}$ in the $k$-lattice.

In order to put the density of counting numbers along $\R_+$ and counting numbers along $\mathbf{R}_{k}\subset\R^k_+$ on equal footing, we use the gamma hypothesis for prime powers together with the fundamental theorem of arithmetic. Evidently, it suffices to deduce their ratio for the particular case of counting prime powers. So the task is to determine the prime-power density along $\mathbf{R}_{k}$ relative to the prime-power density along $\R_+$.

We require the probability that the point $\mathrm{r}_k=(\mathrm{r},\mathrm{r}+h_2,\ldots,\mathrm{r}+h_k)\in\mathbf{R}_{k}$  lies on the pair-wise coprime $k$-lattice and is coprime to some prime since this is a necessary condition for $\mathrm{r}_k$ to be a prime power $k$-tuple. Our main tool is a theorem by T\'{o}th \cite{T}:
Let $k,m,u\geq1$ and
\begin{equation}
P_k^{(u)}(m)=\!\!\!\!\!\!\!\!\sum_{\begin{array}{c}
            \scriptstyle{1\leq a_1, \ldots, a_k\leq m} \\
            \scriptstyle{(a_i,a_j)=1, i\neq j} \\
            \scriptstyle{(a_i,u)=1}
          \end{array}}\!\!\!\!\!\!\!\!1
\end{equation}
be the number of $k$-tuples $(a_1,\ldots, a_k)$ on the pair-wise coprime lattice with $1\leq a_i\leq m$ and $(a_i,u)=1$ for all $i\in\{1,\ldots, k\}$.
\begin{theorem}\emph{(\cite{T})}
For a fixed $k\geq1$, we have uniformly for $m,u\geq1$,
\begin{equation}
P_k^{(u)}(m)=A_k f_k(u) m^k+O\left(\theta(u)m^{k-1}\log^{k-1}(m)\right)
\end{equation}
where\begin{eqnarray*}
       A_k &=& \prod_{p'}\left(1-\frac{1}{p'}\right)^{k-1}\left(1+\frac{k-1}{p'}\right) \\
       f_k(u) &=&\prod_{p'|u}\left(1-\frac{k}{p'+k-1}\right)
     \end{eqnarray*}
and $\theta(u)$ is the number of squarefree divisors of $u$.
\end{theorem}
Restrict $P_k^{(u)}(m)$ to the ray $\mathbf{R}_{k}$ and choose $m>\mathrm{r}+h_k$. Then the density of points along $\mathbf{R}_{k}$ that are coprime to a given prime (or prime power) $p$ is
\begin{equation}
D_k^{(p)}(m):=P_k^{(p)}(m)/m^k=A_k f_k(p)+O\left(\log^{k-1}(m)/m\right)\;.
\end{equation}
Of course $m$ is automatically coprime to all primes $p> m$. Consequently, the density of prime powers along $\mathbf{R}_{k}$ is
\begin{equation}
\lim_{m\rightarrow\infty}\prod_{p\leq m} D_k^{(p)}(m)=\prod_p\left(1-\frac{1}{p}\right)^{k}\;.
\end{equation}
In other words,  the prime-power measure along $\mathbf{R}_{k}\subset\R_+^k$ relative to the prime-power measure on $\R_+$ is
\begin{equation}
\prod_p\left(1-\frac{1}{p}\right)^{k}\;d\mathrm{r} = dr
\end{equation}

Meanwhile, to determine the contribution from the `smoothed' $\Lambda_{(k)}(n)$ to the prime-power measure up to a cut-off $x$, follow the standard argument. Let $\nu_p(\mathcal{H}_k)$ denote the number of distinct residue classes mod $p$ occupied by  $\mathcal{H}_k$.  There are $p-\nu_p(\mathcal{H}_k)$ residue classes mod $p$ that $n$ can occupy that guarantee $n^k_{(k)}$ is not divisible by $p^k_{(k)}< n^k_{(k)}$. So the probability that $\Lambda_{(k)}(n)$ does not vanish given that $n$ is a prime power is
 \begin{equation}
 \prod_{p_{(k)}<n_{(k)}} \left(1-\frac{\nu_p(\mathcal{H}_k)}{p}\right)\;.
\end{equation}

It follows that the prime-power probability measure on $\mathbf{R}_k$ up to some cut-off $x$ is
\begin{equation}\label{prime constant}
\prod_{p_{(k)}< x_{(k)}} \left(1-\frac{\nu_p(\mathcal{H}_k)}{p}\right)\,d\mathrm{r}=\prod_{p< x} \left(1-\frac{\nu_p(\mathcal{H}_k)}{p}\right)\prod_p\left(1-\frac{1}{p}\right)^{-k}\,dr
=:C_{(k)}(x) \,dr
\end{equation}
where $dr$ is the integrator associated with the normalized
Haar measure on $\R_+$. Recall that $\nu_p(\mathcal{H}_k)=k$ as soon as $p>h_k$. Finally, define $C_{(k)}$ to be the asymptote
\begin{equation}
C_{(k)}:=\lim_{x\rightarrow\infty}C_{(k)}(x) =\prod_p\left(1-\frac{\nu_p(\mathcal{H}_k)}{p}\right)\left(1-\frac{1}{p}\right)^{-k}\;.
\end{equation}

\section{Comparing gamma and Hardy-Littlewood}
Let's compare the counts of prime doubles between Conjecture \ref{average primes} and Hardy-Littlewood. First, we note that
\begin{equation}
\frac{1}{\log(r)\log(r+h_{2i})}\sim\frac{1}{\log(r)^2}-\frac{h_{2i}}{r\log(r)^3}
+O\left(\frac{h_{2i}^2}{r^2\log(r)^3}\right)\;.
\end{equation}
For cut-off $x$ and off-set $h_{2i}$,  we don't expect much difference between the two associated integrals when $x\gg h_{2i}$. But for the other way around $h_{2i}> x$, there may be. We employ \emph{Mathematica} to explore some numbers.

The table below contains the exact number of prime doubles $(n,n+h_{2i})$ such that $n\in(2,x)\subset\Z_+$ for several values of $h_{2i}$.
\begin{table}[H]
\centering
\begin{tabular}{l||c|c|c|c|c|c|c|}
  $x\,\diagdown \;h_{2i}$& $10^1$ & $10^2$ & $10^3$ & $10^4$ & $10^5$ & $10^6$ & $10^7$ \\\hline\hline
  $10^1$ & 1 & 2 & 0 & 1 & 1 & 1 & 0\\\hline
  $10^2$ & 2 & 9 & 5 & 5 & 3 & 2 & 2\\\hline
  $10^3$ & 16 & 49 & 37 & 34 & 23 & 20 & 17\\\hline
  $10^4$ & 111 & 260 & 253 & 224 & 186 & 163 & 142\\\hline
  $10^5$ & 859 & 1615 & 1631 & 1556 & 1431 & 1219 & 1050\\\hline
$10^6$ & 6707 & 10906 & 10993 & 10798 & 10629 & 9766 & 8592\\\hline
 $10^7$ & 56334 & 78248 & 78265 & 77850 & 77680 & 76212 & 71247\\\hline
\end{tabular}
\caption{Exact number of prime doubles for the indicated cut-off $x$ and off-set $h_{2i}$.}
\end{table}

For the Hardy-Littlewood estimate, we have for $h_{2i}=10^1$
\begin{equation}
C_{(2)}=\frac{1-\frac{1}{2}}{\left(1-\frac{1}{2}\right)^2}\cdot\frac{1-\frac{2}{3}}{\left(1-\frac{1}{3}\right)^2}
\cdot\frac{1-\frac{1}{5}}{\left(1-\frac{1}{5}\right)^2}
\cdot\frac{1-\frac{2}{7}}{\left(1-\frac{1}{7}\right)^2}
\cdot\prod_{7<p}\frac{1-\frac{2}{7}}{\left(1-\frac{1}{7}\right)^2}\approx 1.76\;.
\end{equation}
All subsequent values of $h_{2i}$ produce the same leading terms in the product since they are all powers of $10$. So $C_{(2)}\approx 1.76$.
Using
\begin{equation}
C_{(2)}\int_2^x\frac{dr}{\log(r)^2}\;,
\end{equation}
the percentage deviation of the Hardy-Littlewood estimate from the exact count
is tabulated below:
\begin{table}[H]
\centering
\begin{tabular}{l||c|c|c|c|c|c|c|}
  $x\,\diagdown \;h_{2i}$& $10^1$ & $10^2$ & $10^3$ & $10^4$ & $10^5$ & $10^6$ & $10^7$ \\\hline\hline
  $10^1$ & 222 & 222 & - & 545 & 545 & 545 & -\\\hline
  $10^2$ & 64 & 101 & 261 & 261 & 502 & 802 & 802\\\hline
  $10^3$ & 19.7 & 24.6 & 65 & 79.6 & 165 & 205 & 259\\\hline
  $10^4$ & 5.78 & 9.85 & 12.9 & 27.5 & 53.6 & 75.2 & 101\\\hline
  $10^5$ & 2.52 & 3.1 & 2.1 & 7 & 16.3 & 36.6 & 58.6\\\hline
$10^6$ & .58 & .84 & .04 & 1.85 & 3.47 & 12.6 & 28\\\hline
 $10^7$ & .16 & .12 & .09 & .63 & .85 & 2.8 & 9.96\\\hline
\end{tabular}
\caption{Percentage deviation between exact and Hardy-Littlewood estimates of prime doubles for the indicated cut-off $x$ and off-set $h_{2i}$.}
\end{table}
Since the Hardy-Littlewood estimate is asymptotic, it is not surprising that percentages are fairly high for lower cut-offs. But notice the general trend of increasing deviation across rows as the ratio $h_{2i}/x$ increases.

Now, for the estimate from Conjecture \ref{average primes} we need
\begin{equation}
C_{(2)}(x)=\frac{1-\frac{1}{2}}{\left(1-\frac{1}{2}\right)^2}\cdot\frac{1-\frac{2}{3}}{\left(1-\frac{1}{3}\right)^2}
\cdot\frac{1-\frac{1}{5}}{\left(1-\frac{1}{5}\right)^2}
\cdot\frac{1-\frac{2}{7}}{\left(1-\frac{1}{7}\right)^2}\cdot\prod_{7<p\leq x}\frac{1-\frac{2}{7}}{\left(1-\frac{1}{7}\right)^2}\;.
\end{equation}
Evidently this pre-factor ranges between $1.76\leq C_{(2)}(x)\leq 1.83$. Using
\begin{equation}
\sum_{m=1}^\infty\frac{\mu(m)}{m}C_{(2)}(x^{1/m}){\int_2}^{\,x^{1/m}}\frac{dr}{\log(r)\log(r+h_{2i})}\;,
\end{equation}
the percentage deviation of the gamma conjecture estimate from the exact count is tabulated below:
\begin{table}[H]
\centering
\begin{tabular}{l||c|c|c|c|c|c|c|}
  $x\,\diagdown \;h_{2i}$& $10^1$ & $10^2$ & $10^3$ & $10^4$ & $10^5$ & $10^6$ & $10^7$ \\\hline\hline
  $10^1$ & 57.2 & 8.2 & - & 7.15 & 25.7 & 38.1 & -\\\hline
  $10^2$ & 7.28 & 1.06 & 30 & 1.88 & 30.9 & 63.6 & 40.25\\\hline
  $10^3$ & 1.84 & .95 & 10.24 & 5.83 & 11.87 & 7.24 & 8.15\\\hline
  $10^4$ & .76 & 3.39 & 1.33 & .75 & .47 & 4.12 & 5.63\\\hline
  $10^5$ & 1.19 & 1.58 & .29 & .92 & .7 & .11 & .3\\\hline
$10^6$ & .23 & .47 & .46 & .68 & .48 & .2 & .46\\\hline
 $10^7$ & .07 & .02 & .02 & .41 & .12 & .2 & .4\\\hline
\end{tabular}
\caption{Percentage deviation between exact and gamma conjecture estimates of prime doubles for the indicated cut-off $x$ and off-set $h_{2i}$.}
\end{table}
At least in these parameter ranges the gamma estimates are superior, but there is no reason not to expect similar comparisons throughout parameter space and for all prime $k$-tuples.

\end{document}